\pageno=1 

\parindent=0cm
\magnification 1200

%

\tolerance 40000000
=%
cmr7
=%
cmbx7
=%
cmbx12

\def \o {{\omega}}

\def \R {{\bf R}}









\def \d {\partial}
\def \del {\partial}

\def \GRAD {\nabla\!_x}
\def \DIV {\nabla\!_x \! \cdot \!}
\def \LAP {\Delta_x}

\def \eq { \eqno}
\def \eps {{\epsilon}}

\def \O {\Omega}
\def \O {{\Omega}}

\def \frac#1#2{{\textstyle{#1 \over #2}}}

\def  \eq  {\eqno}


\centerline {\bf Weak Convergence and Deterministic}
\bigskip
\centerline {\bf Approach to Turbulent Diffusion.}
\bigskip
\bigskip

\bigskip
\centerline { Claude BARDOS
\footnote{${}^1$ }
{ {\rm
Two authors C. Bardos and S. Kamvissis acknowledge the support of the TMR
Asymptotic in Kinetic Theory for this contribution. }}, Jean-Michel
GHIDAGLIA}
\centerline {\it Centre de Math\'ematiques et de Leurs Applications }
\centerline {\it Ecole Normale Sup\'erieure de Cachan et CNRS UMR 8536 }
\centerline {\it 94235 CACHAN  CEDEX  FRANCE}
\centerline {$\{bardos,jmg\}@cmla.ens-cachan.fr$}
\bigskip
\centerline {AND}
\bigskip
\centerline { Spyridon KAMVISSIS${\,}^1$ }
\centerline {\it UNIVERSITY OF AIX-MARSEILLES, FRANCE and MSRI, BERKELEY } 
\centerline {$kamvissi@msri.org$}

\bigskip
\bigskip
\rightline {\it Dedicated to Walter Strauss }
\bigskip
\bigskip
\bigskip
 {\bf Abstract.}
\smallskip\noindent

\bigskip
The purpose of this contribution is to show that some of the basic ideas
of turbulence can be addressed
in a deterministic setting instead of introducing random realizations of
the fluid. 
Weak limits of oscillating sequences of solutions are considered and
along the same line the Wigner
transform replaces the Kolmogorov definition of the spectra of
turbulence. One of the main issue is to show
that, at least in some cases, this weak limit is the solution of an
equation with an extra diffusion (the
name turbulent diffusion appears naturally). In particular for a weak
limit of solutions of the 
incompressible Euler equation (which is time reversible)  such process
would lead to the appearance of
irreversibility.  
In the absence of   proofs, following a program initiated by P. Lax [L],
the 
diffusive property of the limit is analyzed, with the tools of Lax and
Levermore [LL] or Jin Levermore and 
Mc Laughlin [JLM], on   the zero dispersion limit of the
Korteweg-deVries equation and of
the Non Linear Schrodinger equation. The three authors are extremely
happy to have the opportunity to publish
this contribution in a volume dedicated to Walter Strauss as a mark of
friendship and admiration for his
achievement. They hope that this paper concerned with non linear fluid
mechanics, non linear instabilities 
and inverse scattering, will find its place in the different domains
that have interested Walter.
\bigskip
\vfill
\eject 
{\bf 1. Introduction.}
\bigskip
Two types of objects appear frequently in the theory of turbulence for a
fluid defined in an open set $\O
\subset \R^d,\, d=1, 2 \hbox{ or } 3$~:  Models of Turbulence and 
Turbulent
spectra, the first ones is used in most of practical numerical
simulations and the second  
proposed as a tool for the understanding of the phenomena.

One of the most classical  models of turbulence is the so called $k-\eps$
model which reads~:
$$
\eqalign{ \d_t U +\GRAD(U\otimes U) +\GRAD P -\nu \LAP U-
c_\nu\GRAD [{{k^2}\over\eps}(\GRAD U+ ^t\GRAD U)]&=0,\cr
\DIV U&=0\cr}
\eq(1)
$$

and
$$
\eqalign{\d_tk + U\GRAD k -{{c_\nu k^2}\over {2 \eps}}|\GRAD U +(\GRAD
U)^t|^2-
\DIV[ c_\nu{{k^2}\over \eps}\GRAD k]+\eps & =0\,,\cr
\d_t\eps + U\GRAD \eps-{{c_1 k}\over {2 }}|\GRAD U +(\GRAD U)^t|^2-
\DIV[ c_3{{k^2}\over \eps}\GRAD \eps]+c_2{{\eps^2}\over k}& 
=0\,.\cr}\eq(2)
$$
In the above system $c_i$ denote several constants usually given by
experimentation. Equation (1)
is the standard incompressible Navier Stokes equation modified by the
introduction of a "turbulent viscous" term~:
$$
 \nu \GRAD u+  c_\nu[{{k^2}\over\eps}(\GRAD U+ ^t\GRAD U)]\,;\,\, \nu_T 
= \nu +c_\nu{{k^2}\over\eps} \eq (3)
$$   
where the space-time dependent scalar quantities $k$ and $\eps$
are defined by the 
coupled system (2) and interpreted in term of local energy fluctuation
and local enstrophy.

Intuitively the formula encompasses the idea that the  
fluid by its self interaction
produces some averaging effect and therefore that the average $U$ is a
smooth quantity solution of 
an  equation with a viscosity greater than the initial one. Therefore
one conjectures and for some
variant of the system (1) (2) (cf.  [MP] and [LM]) proves the relation~:
$$
\nu_T \ge \nu\,. \eq(4)
$$
This relation becomes  important when $\nu$ is very small. This
situation corresponds to large Reynolds
numbers. In some sense the idea is to have for the average an equation
which would become valid when, due to
the  complexity, the initial equation cannot be computed. 

Such effect becomes crucial when $\nu=0$. All what
remains
is the turbulent aspect of the viscosity. It has to be positive
otherwise the turbulent model is an ill
posed problem (in the sense of Hadamard) and the solution cannot be
computed. Furthermore it gives one of
the many examples where the limit of a reversible system (the
incompressible Euler equation) becomes
irreversible. 

There is up to now no mathematical analysis  of the range of validity of
such formulas.
The "phenomenological" proof usually given introduces some randomness in
the description of the 
motion of the fluid.
The randomness plays an even more crucial role in the definition of the
spectra of turbulence
as given below. A random family of realizations $u(x,t,\o),$ is said to
be homogeneous whenever the tensor
$$
\langle u(x+r,t)\otimes u(x,t) \rangle
$$
is independent of $x$ and $t$ and it is said to be isotropic whenever it
depends only on $|r|$.
The homogeneity hypothesis implies that the tensor 
$$
\widehat R (k)=\int u (x+r,t,\o)\otimes u (x,t,\o)e^{-ir\cdot k}d\o dr
\eq(5)
$$
is independent of $x$ and $t$. The isotropy hypothesis implies that it
is given by the formula~:
$$
\widehat R (k,t)= {{E(|k|)}\over {4\pi |k|^2}} (I-{{k\otimes
k}\over{|k|^2}})\,. \eq(6)
$$
The expressions (5) and (6) are starting points for most of the
phenomenological theory of turbulence
according to Kolmogorov [Ko] and Kraichnan [Kr].

\bigskip

{\bf 2 Weak convergence and Wigner Measures.}
\bigskip
As said in the introduction one of the purpose of this contribution is
to show that the recourse to
randomness is not compulsory  for the above construction. At variance,
as already observed by Lax [L], one
could consider the weak limit of a deterministic sequence of oscillatory
solutions of the fluid mechanics
equations and use recent results concerning defect measures and Wigner
transform (cf.[Ge], [LP] and [Ta]).

The first observation is that with the homogeneity hypothesis, the 
formula (5) can also be written as~:
$$
\widehat R (k,t)=\int u (x+{r\over 2},t,\o)\otimes u (x-{r\over
2},t,\o)e^{-ir\cdot k}d\o dr \eq(7)\,,
$$
and therefore the right hand side turns out to be the Wigner transform
of $u$. As above it is
a symmetric positive tensor. However  the  progress is that this
tensor is can be written in 
term of a local quantity 
$\widehat R (x,t,k,\o)$~: 
$$
\eqalign{ \widehat R (x,t,k)=&<\widehat R (x,t,k,\o)>\cr
\widehat R (x,t,k,\o)=&\int u (x+{r\over 2},t,\o)\otimes u (x-{r\over
2},t,\o)e^{-ir\cdot k}  dr \cr}\eq(8)
$$ 
and  homogeneity hypothesis can be relaxed. 

In fact this idea already appeared as a basic
ingredient
of the contribution of  D.C. Besnard, F.H. Harlow,  R.M. Rauenzahn and
C.
Zemach [BHRZ]. In spite of the fact that it is a very natural approach, 
to
the
best of our knowledge, it has not been used elsewhere in turbulence
theory. Furthermore the Wigner
transform, in space, or (after time localization) in time is well
defined for any  solution $u$ or sequence
of solutions $u_n$ and one has the standard formulas~:
$$
\eqalign{
&\widehat R_n (x,t,k)=\int u_n (x+{r\over 2},t)\otimes u_n (x-{r\over
2},t,\o)e^{-ir\cdot k}  dr\,, \cr
&u_n (x,t)\otimes u_n (x,t)= ({1\over { 2\pi}})^d\int \widehat R_n
(x,t,k)dk\,, \cr
&\widehat R_n^\theta (x,t,\tau)=\int \theta (s) u_n (x,t+{s\over
2})\otimes u_n (x,t-{s\over
2},\o)e^{-is\tau}  ds,\,\, \theta \in {\cal D} (\R_+), \, \theta(0)
=1\,, \cr
&u_n (x,t)\otimes u_n (x,t)=({1\over { 2 \pi}})\int \widehat R_n^\theta
(x,t,\tau)d\tau\,.\cr} \eq(9)
$$
 The two last formulas of (9) are   important because they involve
instead of space correlations,
time correlations which are the quantities  more commonly involved in
practical experiments.

The only physical {\it a priori} estimate, uniform with respect to the
Reynolds number, is the energy estimate.
$$
\frac12 \int_\O |u_n(x,t)|^2dx \le \frac12 \int_\O |u_n(0,t)|^2dx \le
C\,. \eq(10)
$$
This observation is valid in particular for the $3d$ Euler or Navier
Stokes equation as proven by Di Perna
and Lions (cf [Li] section 4.3). 
It is also valid even in $2d$ for the Navier Stokes equation with
viscosity going to zero, when the natural
viscous boundary condition $u=0$ on $\d \O$ is assumed. In particular 
Grenier has constructed solutions to
the Navier Stokes equation in the half plane $\O =\R_{y_+}\times \R_x$
with vorticity blowing up in the
$L^\infty$ norm when the viscosity goes to zero([Gre] theorem 2.1.). As
observed by Grenier this is an 
instability phenomenon for the Prandtl layer which is of the same type as
the non linear instability of the
Euler equation  proven by Friedlander Strauss and Vishik [FSV].

With the estimate (10) one 
concludes that, up to the extraction of a subsequence,
$u_n$ converges in $weak^*L^\infty (\R_{t_+}, L^2(\O))$ to a limit $u$.
However due to the above
considerations, in many cases, one will have~:
$$
\lim_{n\rightarrow \infty} u_n (x,t)\otimes u_n (x,t) = u (x,t)\otimes
u  (x,t) + R_{turb}(x,t) \eq(11)
$$
with $R_{turb}(x,t) \not=0$. In a follow up of ideas of Peter Lax [L]
one could consider that the appearance
of the tensor
$R_{turb}$  plays the role of the Reynolds stress tensor as the emergence
of turbulence in a
deterministic approach. In
fact it is the Defect Measure of the sequence $u_n$. The  tensor
$R_{turb}$ is symmetric positive definite and
one has~:
$$
R_{turb}(x,t) = \lim_{n\rightarrow \infty}(u_n-u)\otimes(u_n-u)\,,
\eq(12)
$$
or, with the Wigner transform of $u_n-u$,
$$
\widehat R_{turb}(x,t,k) = \lim_{n\rightarrow \infty} 
\int_{\R^d} (u_n-u)(x+{r\over 2},t)\otimes(u_n-u)(x-{r\over
2},t)e^{-ir\cdot k}dk \eq(13)
$$
and
$$
\lim_{n\rightarrow \infty} u_n (x,t)\otimes u_n (x,t) = u (x,t)\otimes
u  (x,t) + ({1\over {2 \pi}})^d\int_{\R^d}\widehat
R_{turb}(x,t,k) dk\,. \eq(14)
$$
The formula (14) displays the natural link, for our purpose, between the
defect measure and the Wigner
transform. Observe that the left hand side of (13) is the natural local
and deterministic ``avatar" of the
Kolmogorv spectra for random turbulence and it is natural to conjecture
that it will inherit the basic
properties of isotropy and scale law for the dependence in $k$.
The isotropy hypothesis is made plausible $3d$ by the following remarks~:

A necessary condition for $R_{turb}(x,t)$ to be non zero is that $\hbox{
curl } \,u_n$ becomes unbounded in
the neighborhood of $(x,t)$ then~:

(i) for non zero viscosity it  has been shown by Constantine and
Fefferman [CF] that it is   much more
the oscillations in  direction of the vorticity than its size that are
responsible for instabilities in the
fluid,

(ii) exterior constant Coriolis force  stabilizes the fluid when the
Rossby number goes to infinity as
shown by Babin Nicolaenko and Mahalov [BNM] and others. Therefore a
decomposition of the vorticity according
to the formula~:
$$
\hbox{ curl }u_n = \O_n + \tilde \o_n \eq(15)
$$
with $\O_n $ having a constant direction and a modulus going to $\infty$
while $\o_n$ remains bounded
should not be possible.  

\bigskip
In the sequel of this section we consider in the two dimensional case,
sequences of solutions to
the Euler equation with the impermeability condition~:
$$
\eqalign{ \d_t u_n  + \GRAD( u_n\otimes u_n) +\GRAD p_n &=0\,,\cr
\vec n \cdot u =0 \hbox{ on }  \d \O,\,\, \DIV u_n&=0 \,,\cr}
\eqno(16)
$$ 
or  of the Navier Stokes equation in a domain
$\O$ of $\R^2$  with the viscous boundary condition~: 
$$
\eqalign{ \d_t u_n +\GRAD( u_n\otimes u_n)  +\GRAD p_n -\nu_n \LAP
u_n&=0,\cr
u_n=0 \hbox{ on }  \d \O,\,\, \DIV u_n&=0\,.\cr}
\eq(17)
$$
 
In both cases  existence and uniqueness of such solutions are well
established facts with the
hypothesis
$$
u_n(x,0)=u_n^0(x) \in L^2(\O), \DIV u_n^0 =0, \o_n^0=\nabla\times u_n^0
\in L^\infty(\O)\,, \eq(18)
$$
in the first case and with the assumption
$$
u_n(x,0)=u_0(x) \in L^2(\O), \DIV u_0 =0\,,
$$
in the second case ($\nu_n > 0$).

\bigskip
However in these two cases the turbulent Reynolds tensor may be present
in the limit~:

(i) if one considers a sequence $u_n$  of solutions of the Euler
equation with  initial
data uniformly bounded in $L^2(\O)$  but with initial  $\o_n^0$
vorticity unbounded in $L^\infty$,

(ii) in  the second case if, as already mentioned above, one keeps the
initial data fixed but let the
viscosity $\nu_n$ goes to zero. 

In both cases the limit satisfies the equations
$$
\d_t u + \GRAD (u\otimes u) +\GRAD R_{turb}+ \GRAD p = 0, \, \DIV u
=0\,, \eq(19)
$$
with
$$
R_{turb}(x,t)=\lim_{n\rightarrow \infty}
(u_n(x,t)-u(x,t))\otimes (u_n(x,t)- u(x,t))\,.\eq(20)
$$
 Introducing the trace~:
$$
T = {{R_{turb}^{11} +R_{turb}^{22}}\over 2}
$$
(19) is changed into 
$$
\d_t u + \GRAD (u\otimes u) +\GRAD S_{turb}+ \GRAD P = 0, \, \DIV u
=0\,. \eq(21)
$$
with $S_{turb}$ denoting a tracefree tensor and $P=p+T$. The space of
tracefree tensors is of dimension $2$ 
and assuming that
that the limit $u$ is a smooth function, this space has a natural basis
given by the matrix~:
$$
\frac12\big( \nabla u + \nabla^t u\big) = \left(\matrix{{\del_{x_1}
u_1}, & \frac12(\del_{x_2} u_1
+\del_{x_1} u_2)\cr
\frac12(\del_{x_2} u_1 +\del_{x_1}
u_2)&\del_{x_2} u_2\cr}\right)  \eq(22)
$$
and an orthogonal complement
$$
\Phi (u) =\left(\matrix{\frac12(\del_{x_2} u_1 +\del_{x_1}
u_2)&\del_{x_2} u_2 \cr
\del_{x_2} u_2 &-\frac12(\del_{x_2} u_1 +\del_{x_1}
u_2)}\right)\,. \eq(23) 
$$
Therefore there exist  two space-time depending functions
$\nu_{turb}(x,t)$ and $\delta(x,t)$ such that one has~:
$$
S_{turb}= \nu_{turb}\big( \nabla u + \nabla^t u\big) + \delta \Phi(u)
\eq(24)
$$
and the equation (21)
becomes the equation~:
$$
\d_t u + \GRAD( u\otimes u) +
\GRAD\big(\nu_{turb} \big( \nabla u + \nabla^t u\big)\big)  +\GRAD
(\delta \Phi(u))  + \GRAD P = 0, \, \DIV u
=0\,.
\eq(25)
$$
As discussed above, the mechanism of "creation of turbulence" should be
isotropic and this would imply   that
the tensor
$S_{turb}$ is invariant under Galilean transformations and therefore 
proportional to $\big( \nabla u +
\nabla^t u\big)$ reducing (25) to a diffusive type equation~:
$$
\d_t u + \GRAD( u\otimes u) +
\GRAD\big(\nu_{turb} \big( \nabla u + \nabla^t u\big) \big) +\GRAD P =
0, \, \DIV u =0\,.
\eq(26)
$$

Observe however that in any case, assuming that
on the boundary
$u$ is zero whenever the
$\nu_{turb}$ is positive one has the energy estimate~:
$$
\frac12\del_t\int_\O |u(x,t|^2dx + \int_\O \nu_{turb} (x,t) |\GRAD u|^2 dx
=0\,. \eq(27)
$$

A turbulent model would be obtained by coupling   equations (25) with a
system of equations which would
determine the function $\nu_{turb}$. The necessary condition to obtain a
such well posed system is that
the function 
$
\nu_{turb}(x,t)
$
is non  negative, a property which does not follow   for the fact that
the tensor $R_{turb}$ is itself non
negative.

In the case of the Euler equation, with a sequence of 
initial data $u_n^0$ having unbounded vorticity  one would derive an 
irreversible problem as the limit of
reversible equations. Eventually from the formula (27) one deduces the

\bigskip
{ \bf Proposition 1} In the above configurations, 
assume that the sequence of initial data $u_n^0(x)$ converges 
strongly to $u_0(x)$ in $L^2(\O)$, then  for $T>0,$
the following assertions are equivalent~:

(i) the sequence $u_n(x,t)$ converges strongly in $L^2(\O \times
[0,T])$,

(ii)  $\nu_{turb}(x,t)$ is identically
zero on
$\O
\times [0,T]$,

(iii) one has~:

$$
\int_0^T\int_\O \nu_{turb} (x,t) |\GRAD u|^2 dxdt \le 0\,. \eq(28)
$$

{\bf Proof} The only non classical point is the fact that (iii) implies
the strong convergence or 
equivalently that one has~:

$$
\liminf_{ n \rightarrow \infty} \frac12\int_0^T\int_\O |u_n(x,t)|^2dxdt
=
\frac12\int_0^T\int_\O |u(x,t)|^2dxdt\,. \eq(29)
$$
With the classical energy estimate (for any given $n$) and the relation
(27) one deduces the inequalities~:
$$
\eqalign{{T\over 2}\int_\O |u_0(x)|^2dx \ge \liminf_{ n \rightarrow
\infty} \frac12\int_0^T\int_\O
|u_n(x,t)|^2dxdt =
\frac12\int_0^T\int_\O |u(x,t)|^2dxdt\,,&\cr \ge {T\over 2}\int_\O
|u_0(x)|^2dx -\int_0^T\int_\O \nu_{turb}
(x,t) |\GRAD u|^2 dxdt\,,&\cr} \eq(30)
$$
and (29) follows from (28).

\bigskip

{\bf 3  Positivity versus non positivity of the Diffusion 
Coefficient for the small Dispersion limit of  
KDV and NLS flows.}
\bigskip

In the absence of a systematic theory it seems worth while to study,
following the program of   Lax,
the issue of the positivity of the turbulent coefficient on the
dispersive limit of the KDV and NLS
equations using explicit formulas given by the inverse scattering
theory. 

For the KDV flow one considers the problem~:
$$
u^\eps_t -6 u^\eps u^\eps _x + \eps^2 u^\eps_{xxx}=0,\, \hbox { with
initial data }u^\eps(x,0)=u_0(x) \eq(31)
$$
and the for the NLS flow the problem~:
$$
\eqalign{ i\eps u^\eps_t +{{\eps^2}\over 2} u^\eps_{xx} +
(1-|u^\eps|^2)u^\eps=0&,\cr
\hbox { with initial data } u^\eps(x,0)=A(x)exp\big(i{{S(x)}\over
{\eps}})&\,.\cr}\eq(32)
$$
With the introduction of the 
functions~:
$$
\rho^\eps  =|u^\eps|^2-1, \hbox{ and } \mu^\eps = {-{i\eps}\over 2}
(u^\eps\bar{u^\eps_x}
-u^\eps_x\bar{u^\eps})\,, \eq(33)
$$
the NLS equation is equivalent to the system~:
$$
\eqalign{{\rho^\eps}_t + {\mu^\eps}_x =& 0,\cr
{\mu^\eps}_t + ({{{\mu^\eps}^2} \over {\rho^\eps}} + 
{{{{\rho^\eps}^2} \over 2})}_x  = &{{\eps^2} \over 4}
({\rho^\eps}  {{(\log{\rho^\eps})}_{xx})}_x,\cr}\eq(34)
$$
The equations (31) and (32) are time reversible and (34)  is a
reversible perturbation of the the usual
isentropic compressible Euler equation. For
$\eps$ going to zero
  the
functions
$u^\eps$,
$\rho^\eps$ and
$\mu^\eps$ converge weakly and the following notations are introduced~:
$$
\eqalign{
\bar u &= weak \lim_{\eps \rightarrow 0} u^\eps\,,\cr
\bar{ u^2}&=weak \lim_{\eps \rightarrow 0} {(u^\eps)}^2\,,\cr
\bar{ \rho}&=weak \lim_{\eps \rightarrow 0} \rho^\eps\,,\cr
\bar { \mu } &=weak \lim_{\eps \rightarrow 0} \mu^\eps\,,\cr
Q(\rho^\eps,\mu^\eps) &=({{{\mu^\eps}^2} \over {\rho^\eps}} + 
{{{\rho^\eps}^2} \over 2})\,,\cr
\bar Q& = weak \lim_{\eps \rightarrow 0}
Q(\rho^\eps,\mu^\eps)\,.\cr}\eq(35)
$$
and one obtains the equation~:
$$
\del_t \bar u -6 \del_x (   {{{\bar u}^2}\over 2})-6 \del_x(    {{\bar{
u^2}} \over 2}
- {{{\bar u}^2}\over 2})=0 \eq(36)
$$
and the system
$$
\eqalign{\bar\rho_t  + \bar \mu _x =& 0\,,\cr
\del_t{\bar \mu } + \del_x({{\bar \mu ^2} \over {\bar\rho}} + 
{{\bar \rho ^2} \over 2}) +\del_x(\bar Q -Q(\bar \rho,\bar \mu))
 = &0\,.\cr} \eq(37)
$$
In the region where strong convergence occurs one has
$$
(    {{\bar{ u^2}} \over 2}
- {{{\bar u}^2}\over 2})=0 \eq(38)
$$
or
$$
\bar Q -Q(\bar \rho,\bar \mu))
 =  0\,. \eq(39)
$$
As expected in these regions (36) is  (up to a simple change in the $x$
variable ) the Burgers equation and
(37) is the compressible Euler equation for isentropic fluids.

 On the other hand it is known  (cf [LL] and [JLM]) that the
strong convergence does not hold  everywhere. The regions where strong
convergence fails are called the
{\it Whitham regions}. By a simple convexity argument, (observe that
both the functions
$$
u\rightarrow u^2 \hbox { and } (\rho, \,\mu) \rightarrow Q(\rho,\mu)
$$
are convex), one has  in the Whitham region~:
$$
{{\bar{ u^2}} \over 2}
- {{{\bar u}^2}\over 2} > 0 \hbox { and } \bar Q -Q(\bar \rho,\bar \mu))
> 0\,. \eq(40)
$$

To analyze the possibility of the appearance of "turbulent viscosity."
one writes (36) and (37) in the 
following form~:
$$\bar u_t -3((\bar u)^2)_x -\del_x(\nu_{turb}\del_x \bar  u )=0, \hbox
{ with }  \nu_{turb} (x,t)=
{{\bar{u^2}  -(\bar u)^2} \over {\del_x \bar  u}} \eq(41) 
$$
and
$$
\eqalign{\bar\rho_t + \bar\mu_x = 0\,,&\cr
{\bar \mu }_t + ({{\bar \mu^2} \over {\bar \rho}} + 
{{\bar \rho^2} \over 2})_x - \del_x(\nu_{turb} \bar \mu_x)=0\,,&\cr
\hbox { with } \nu_{turb}(x,t) = {{\bar Q -Q(\bar \rho,\bar \mu))}\over
{\del_x\bar \mu}}\,.&\cr} \eq(42)
$$
The existence of any kind of turbulent model requires that  $\nu_{turb}$
be non negative which is equivalent
here to the property~:
$$
\del_x \bar u (x,t) > 0\,,  \eq(43)
$$
in the Witham region for the KDV dispersive limit and~:
$$
{{\bar Q -Q(\bar \rho,\bar \mu))}\over {\del_x \bar \mu}}\ge 0\,, 
\eq(44)
$$
in the Witham region for the dispersive NLS limit. This is also
equivalent to
$$
\del_x \bar
\mu>0\,. \eq(45)
$$

 Since diffusion properties may appear on a larger time scale it is
natural to explore
the properties (43) and (45) for large time. 

Such program is done below using the tools of the inverse scattering
following [LL] and the conclusion will
be   the fact that such properties are satisfied depends on the initial
data.
The starting point are the following theorems~:

{\bf Theorem 2.} ([LL])  Let $u(x,t; \epsilon)$ solve
$$
u_t -6uu_x +\epsilon^2 u_{xxx}=0, \eq(46)
$$
with initial data $u(x,0; \epsilon) = u_0(x)$ belonging to the Schwartz
class, strictly negative, and with only one minimum point $x_0$, 
at  which $u_0(x_0) =-1$.
Let $x_{\pm} (\eta)$ be defined for $0 < \eta <1 $ by 
$$u_0 (x_{\pm} (\eta)) = -\eta^2\,\hbox{  and }x_- < x_0 < x_+\,.\eq(47)
$$
Define the function
$$
\phi (\eta) = \int^{x_+(\eta)}_{x_-(\eta)}
{{ \eta } \over {(-u_0(x)-\eta^2)^{1/2}}} dx\,, \eq(48)
$$
for $ 0 < \eta < 1$. Then,

(i) the weak limit 
$$\bar u(x,t) = \lim_{\epsilon \to 0} u(x,t ; \epsilon)
$$
exists.

(ii) As $t$ goes to infinity, for $x$ such $\delta < x/t < 4- \delta$,
with $\delta$ is any 
given small
positive constant one has
$$
\bar u(x,t) = -{1 \over {4 \pi t}} \phi (( {x \over {4t}})^{1/2} ) +
o(1/t)\,. \eq(49)
$$

(iii) As $t$ goes to infinity, for $x/t <0$ or $x/t >4$ one has 
$$
\bar u = O(t^{-2})\,.\eq(50)
$$

This theorem  is stated and proved in the third paper of Lax and
Levermore [LL] pages 810-815. 
Furthermore, with some conjecture on the uniform effect of the remote
part of the initial data on  the
solution of KDV equation the authors adapt their asymptotic analysis to
the initial data~:
$$
u_0(x)=  \cases{ -1 &if $x<0$;\cr 0 &if $x>0.$ \cr}\eq(51)
$$
which correspond to the shock profile for the Burgers equation. They
obtain for the weak limit 
the formula~:
$$
\eqalign{ \bar u(x,t)=& -1 \hbox { for }x<-6t\,,\cr 
\bar u(x,t)=& s({{x} \over {t}})\hbox { for } -6t<x<4t\,,\cr
\bar u(x,t) =& \hbox { for } 4t <x\,,\cr}\eq(52)
$$
where the function $\xi \mapsto s(\xi)$ can be computed in term of
complete elliptic integrals. Explicit 
numerical computations done in [LL] on the formula for $s(\xi)$ indicate
that this is an increasing
function on the interval
$[-6,4].$

The theorem 2 has a  counterpart for the NLS dispersive limit.

{ \bf Theorem 3.}  Let $u(x,t;\eps)$ solve the NLS flow~:
$$
i\eps u_t (x,t;\eps) + {\eps^2 \over 2} u_{xx}(x,t;\eps) +
(1-|u(x,t;\eps)|^2) u(x,t;\eps) = 0\,, \eq(53)
$$
$$
u^\eps(x,0)=|u_0(x)|^2 exp\big(i{{S(x)}\over {\eps}})\,,\eq(54)
$$
with
$$\rho_0(x) = 1+ |u_0(x)|^2 \,,\eq(55)
$$ 
and 
$$
\mu_0(x)= \del_x S(x)\,, \eq(56)
$$
belonging to  the Schwartz class. Let
also assume, for simplicity, 
that the initial data are  "single well" in the  following sense
(cf [JLM]). Introduce the functions  
$$
r_{\pm} (x) =  \frac12\del_x S  \pm A(x)\,, \eq(57)
$$ 
and assume that  $r_-$   has only one
maximum   $\lambda_{ max }$ while  $r_+$ has only one minimum 
$\lambda_{min}$ with the relation~:
$$
-1 \le r_-(x)\le  \lambda_{max} < \lambda_{min}\le r_+(x) \le 1 \,.
\eq(58)
$$ 
Define the numbers $x^\pm(\lambda)$ according to the formula~:
$$
\eqalign{ \hbox { for } -1 \le \lambda \le  \lambda_{max},r_-
(x_{\pm}(\lambda)) = \lambda,x_- <x_+\,,\cr
\hbox { for } \lambda_{min} \le \lambda \le 1,r_+ (x_{\pm}(\lambda)) =
\lambda,x_- <x_+\,.\cr}
 $$
Then, the weak limit

$$\bar \rho (x,t) = lim_{\eps \to 0} |u(x,t;\eps)|^2+1 $$ 
exists. Furthermore in  the Whitham region 
$${x\over t} \in (-1,\lambda_{max}) \cup (\lambda_{min}, 1)\,, \eq(59)
$$ 
as $t \to \infty,$
one has~:
$$\eqalign{
\bar \rho (x,t) = 1 - {4 \over {\pi t}} \phi (x/t) (1-(x/t)^2)^{1/2},
\hbox { where}\cr
\phi (\lambda) = \int_{x_-(\lambda)}^{x_+(\lambda)} {{
\lambda - 1/2 (r_+(s) + r_-(s))} \over { (\lambda-r_+(s))^{1/2} (\lambda
- 
r_-(s))^{1/2}}} ds,\cr}\eq(60)
$$
and elsewhere one has, $\bar \rho \simeq 1$.

{\bf Proof}~: The proof is given in [K] and for sake of completeness the
main steps are recalled here.
The existence of
the weak limit is proved in [JLM]. Following a method  suggested in 
[LL]
one begins with the multisolitons formula for
fixed $\eps$ and then let $\eps$ go to zero.

For fixed $\eps$, the long time behavior of $|u|^2$ is as follows [FT,
pp.168-176].
In the solitonless region $|x/t| >1$ or $\lambda_{max} <x/t<
\lambda_{min}$,
one has 
$$|u|^2 =  1 + O (t^{-1/2})\,. \eq(61)
$$ 
As $t \to \infty$. In the Whitham
region, the solution is a multisolitons
solution~:
$$
\eqalign{ |u(x,t;\eps)|^2 \sim 1 - \Sigma_{n=1}^{N(\eps)} s (x-\eta_n t-
x_n, \eta_n),\hbox{ where }\cr
s(x, \eta) =  {{1-\eta^2}\over{cosh^2( (1-\eta^2)^{1/2} {x \over
{2\eps}})}}\,,
\cr}\eq(62)
$$
with exponentially small error. The $\eta_n$'s are the associated 
eigenvalues (of the underlying Dirac operator)
and the $x_n$'s are some phase constants of no importance.

The width of each soliton is 
$$O ({\eps \over{ (1-\eta^2)^{1/2}}})\,. \eq(63)
$$  
By Weyl's law
for the distribution of eigenvalues in $(-1, \lambda_{max}) \cup
(\lambda_{min}, 
1)$ as $\eps \to 0$,
$$
\eta_{n+1} -\eta_n = {{\pi \eps}\over {\phi (\bar \eta_n)}}\,,\eq(64)
$$
where $\bar \eta_n \in (\eta_n, \eta_{n+1}).$

Peaks of solitons are located at $\eta_n t$.
As $t \to \infty$, they are  separated by 
$${{\pi \eps t} \over { \phi (\eta_n) }}\eq(65)
$$  
so for large
$t$ they are well separated. 

The wave number $\eta$ of the soliton that peaks at 
$x$ at time $t$ is
$\eta =x/t$, if $t$ is large and either $-1  < x/t < \lambda_{max}$ or
$\lambda_{min} < x/t < 1.$ 

Therefore the density of the solitons is
$$
{{\phi(x/t)}\over {\pi \eps t}}\,. \eq(66)
$$
The area between  a soliton and the line $u=1$ is
$$
4\eps (1-\eta^2)^{1/2}  \sim 4\eps (1-(x/t)^2)^{1/2}\,,\eq(67)
$$
so the asymptotic area density which is given by $1- \bar \rho$ is  the
product of (66) and (67)~:
$$
{{ 4 \phi (x/t)}\over \pi } (1-(x/t)^2)^{1/2}\,.\eq(68)
$$
Hence,
the asymptotic formula for the weak limit $\bar \rho$ follows.
\bigskip
>From the above statement several observations can be made concerning the
appearance of a positive 
turbulent viscosity in the limit equation satisfied by $\bar u$ for the
KDV dispersive limit and by $(\bar
\rho,
\bar
\mu)$ for the NLS dispersive limit. 
Such positivity would be related to the appearance of irreversibility in
a weak limit
of reversible models.   

As said above the construction of [LL] (section 7 page 817) shows that a
shock profile  as initial data
produces in the limit a smooth solution with
a "turbulent viscosity". 

Following Theorem 2, one considers 
initial data $u(x,0; \epsilon) = u_0(x)$ belonging to the Schwartz
class, strictly negative, and with only one minimum point $x_0$, 
at  which $u_0(x_0) =-1$. For large $t$ the solution is asymptotic to
$$
 -{1 \over {4 \pi t}} \phi (( {x \over {4t}})^{1/2} )  \hbox { with }
\phi (\eta) = \int^{x_+(\eta)}_{x_-(\eta)}
{{ \eta } \over {(-u_0(x)-\eta^2)^{1/2}}} dx\,.\eq(69)
$$
Therefore one has~:
$$
\del_x \bar u \simeq  -{1\over{  \pi
(4t)^\frac32}}\phi'({x\over{(4t)^\frac12}})\,. \eq(70)
$$
And the "turbulent diffusion" hypothesis requires that 
$$
\phi'(\eta)<0 \,, \forall \eta \in [0,1]\,.\eq(71)
$$
With $u_0(x) = -e^{-|x|^\beta} $ one has~:
$$
\phi (\eta) = 2\int^{\big(2\log {1\over \eta}\big)^{1\over \beta}}_{0}
{{ \eta } \over {(-u_0(x)-\eta^2)^{1/2}}} dx\,. \eq(72)
$$
On the table 1 the values of $\phi(\eta)$ 
and for the following $\beta$ exponents ~: $\beta = 1, \, 3/2, 2, \hbox
{ and } 4$ and
with $9$ steps  $ \eta = k 10^{-1}, 1\le k \le 9 $

In the first case $\phi$ is decreasing, in the second case its variation
changes, and then
$\phi$ is increasing in the
two last cases.

For the NLS flow  one considers initial data satisfying the hypothesis
of Theorem 3 and observes that
the existence for large time of a        diffusive regime would be given
by
$
\del_x \bar \mu \ge 0 
$
in the Witham region and   with the conservation law~:
$$\bar\rho_t + \bar\mu_x = 0\,,$$
this is equivalent to
 $\bar \rho_t \le 0
$ in the same region. With Theorem 3 this
condition is equivalent, for large $t$ to the relation
$$
\del_\lambda \big(\lambda \phi(\lambda)(1-\lambda^2)^\frac12 \big)\ge 0,
\forall \lambda \in [-1,+1]\,.
\eq(73)
$$
with $\phi(\lambda)$ given by
$$
\phi (\lambda) = \int_{x_-(\lambda)}^{x_+(\lambda)} {{
\lambda - 1/2 (r_+(s) + r_-(s))} \over { (\lambda-r_+(s))^{1/2} (\lambda
- 
r_-(s))^{1/2}}} ds\,.\eq(74)
$$
Solution with initial  data   having zero momentum and a symmetric
density with an fractional 
exponential rate of convergence at infinity  are
analyzed~:
$$
\rho_0(x) =  1-\frac12 e^{-|x|^\beta}, \,\,\,S(x)=0\,. \eq(75)
$$
The fact that the momentum is zero gives
$r_+(\lambda)+ r_-(\lambda) =0$. The $x$ symmetry of the density remains
true for all time and all 
$\eps$ therefore it is enough to consider the behavior of
$$
\big(\lambda \phi(\lambda)(1-\lambda^2)^\frac12 \big)\phi(\lambda),
\lambda={x\over t}\,,
$$
for 
$$0<\lambda_{min} \le \lambda={x\over t} \le 1\,. \eq(76)
$$
One has~:
$$
\lambda \phi(\lambda)(1-\lambda^2)^\frac12 \big)\phi(\lambda) = 2
\int_0^{\big(-\log
2(1-\lambda)\big)^{1\over \beta}}{{\lambda^2
\sqrt{(1-\lambda^2)}}\over{\lambda^2-(1-\frac12
e^{-|x|^\beta})^2}}dx\,. \eq(77)
$$
The numerical computation given on the table 2 are done for
$\beta=1.5,\,2\,,3$ and $3.5$ with
$\lambda={x/t}$
varying from $\lambda_{min} =0.5$ to  $0.9$ with step $0.1$.
They indicate that the  "turbulent regime" appears for 
$\beta =1.5$, $\beta= 2$, $\beta= 3$ but does not hold for $\beta=3.5$.

\bigskip

$$\vbox{ \settabs 4 \columns \+$\beta=1$&$\beta=
3/2$&$\beta=2$&$\beta=4$\cr
\+5.88251&2.19345&1.75226&.78565\cr
\+5.47775&3.65183&2.07024&1.09422\cr
\+5.06441&3.37626&2.29727&1.36578\cr
\+4.63711&3.09141&2.47717&1.62689\cr
\+4.18879&2.79252&2.62703&1.89375\cr
\+3.70918&2.47278&2.75579&2.18372\cr
\+3.18159&2.12106&2.86860&2.52387\cr 
\+2.57400&1.71600&2.96898&2.97357\cr
\+1.80404&1.20273&3.05934&3.73515\cr}
$$

\centerline{ \bf Table 1 Numerical computation for the dispersive KDV
limit.}

Values of $\phi(\eta)$ are computed for $\beta =1,$  $\beta= 3/2,$ which
appears as a critical case,
$\beta =2$ and $\beta = 4$ with $\eta = 10^{-1} k , \, 1\le k \le 9$
\bigskip
\bigskip

$$\vbox{ \settabs 4 \columns \+$\beta=1.5$&$\beta=
2$&$\beta=3$&$\beta=3.5$\cr
\+0        &0         &0   &.0\cr
\+1.10957&1.29635&1.48400&1.53471\cr
\+1.62180&1.64627&1.63429&1.62240\cr
\+2.16222&1.97843&1.76762&1.70143\cr
\+2.69501&2.21615&1.76986&1.64733\cr}
$$

\centerline{ \bf Table 2 Numerical computation for the dispersive NLS
limit.}

Values of $\lambda\phi(\lambda)$ are computed for $\beta =1.5$  $\beta=
2,$ 
$\beta =3$ and $\beta = 3.5$ which appears as a critical
case with $\lambda = 0.5,\, 0.6,\,0.7,\,0.8,\,0.9.$

\bigskip
{\bf Conclusion.}
\bigskip
In this contribution it has been shown that some of the basic questions
of the statistical
theory of turbulence could be formulated in a deterministic setting with
the introduction 
of sequence of weakly converging solutions. The counterpart of the
turbulent spectra being the
Wigner transform and the turbulent diffusion being related to defect
measures. 
Explicit computations done on integrable classical integrable system
indicate that for these
models it is not always possible to construct a "formal" turbulent
equation. At this point of our analysis
it depends on the behavior of the initial data and in   particular on
the fact that they should not be
too much concentrated (their decay for $|x|$ going to infinity has to be
not too small). 
It is worth while to notice that with the convenient conjectures of [LL]
the shock profile leads always to a
diffusive regime.
\bigskip
\bigskip
\centerline { \bf   BIBLIOGRAPHY}

\bigskip

[BMN1] A. Babin, A. Mahalov and B. Nicolaenko~: Global regularity and
integrability of 3D Euler and Navier-Stokes equations for uniformly
rotating
fluids {\it Asymptotic Analysis} {\bf 15}, No. 2 (1997), 103--150.

[BHRZ] D.C. Besnard, F.H. Harlow,  R.M. Rauenzahn and C.
Zemach~: Spectral transport model for turbulence, {\it Los Alamos
Technical Report} LA-11821-MS.

[CF] P. Constantin and C. Fefferman~:  Direction
of vorticity and the problem of global regularity for the Navier-Stokes
Equations, {\it Indiana University Mathematics Journal } {\bf 42}, No. 3
(1993).

[FT] L.D. Fadeev, L.A. Takhtajan~: Hamiltonian Methods in the Theory of
Solitons,
Springer-Verlag, 1987.

[FSV] S. Friedlander, W. Strauss, M. Vishik~: Non linear instability in
an ideal fluid, {\it
Ann. Inst. H. Poincar\'e Anal. Non lin\'eaire}, {\bf 14} (1987) 187-209.

[Ge] P. G\'erard~: Microlocal Defect Measures, {\it
Comm. PDE}, {\bf 16} (1991) 1761-1794.

[Gre] E. Grenier~: Non d\'erivation des \'equations de Prandlt, 
I-III,
{\it S\'eminaire de l'Ecole Polytechnique}, 17 Mars 1998.

[JLM] S. Jin, C. D. Levermore, and D. W. McLaughlin~: The Semiclassical
Limit of the Defocusing NLS Hierarchy,  {\it Comm.Pure Appl. Math}. 
(1999), to appear.

[K] S. Kamvissis~: Long Time behavior for SemiClassical NLS to appear in
{\it Applied Math. letters}

[Ko] A.N. Kolmogorov~: The local structure of turbulence in
incompressible
viscous fluid for very large Reynolds,{\it C.R. Acad. Sci. URSS} {\bf30}
(1941),
301.

[Kr] R.H. Kraichnan~: Inertial ranges in two-dimensional turbulence, 
{\it
Phys. Fluids} {\bf10} (1967), 1417--1423.

[L] P.D. Lax. The Zero Dispersion Limit, A Deterministic Analogue of
Turbulence,
{\it Comm. in Pure and Appl. Math}, {\bf 54},   (1991), 1047-1056. 

[LL] P.D.Lax. C.D.Levermore~: The Zero Dispersion Limit for the KdV
Equation, 
I-III,
{\it Comm. in Pure and Appl. Math}, {\bf 36}  (1983), 253-290, 571-593,
809-829.

[LM] R. Lewandowsky, B. Mohammadi~: Existence and positivity results for
the $\phi-\theta $ model
and a modified $k-\eps$ model .
{\it Math. Models and Methods in Applied Science},{\bf 3}  (1993),
195-215, 

[Li]  P.L. Lions~:  Mathematical Topics in  Fluid Mechanics Volume 1
Incompressible models, {\it Oxford Lecture Series in Mathematics and
its Applications}, Oxford (1996).

[LP]  P.L. Lions T. Paul~:  Sur les Mesures de Wigner, 
{\it Revista Mat. Iberoamericana}, {\bf 9}, 1993, 553-618.

[MP]   B. Mohammadi and O. Pironneau~:  Analysis of the
K-$\epsilon $ turbulence model,  {\it Research in Applied Math,} No 31,
J.L.Lions and P. Ciarlet (eds.), Masson-Wiley, Paris (1994).

[T] L. Tartar,  H-measures, a new approach for
studying homogeneization, oscillations and concentration effects in
partial differential equations, {\it Proceedings of the Royal Society of
Edinburg,} {\bf 115 -A} (1990), 193--230.

[V] S.Venakides~: Higher Order Lax-Levermore Theory, {\it Comm. in Pure
and Appl. Math}, {\bf 43},  (1990),
335-362.
\bigskip
\bye
\end